\documentclass[preprint,12pt]{elsarticle}

\usepackage{amsopn,amsmath,amssymb,eucal,url,enumerate,amscd}
\usepackage[pagebackref]{hyperref}

 \biboptions{comma,square}
\bibliography{refs}

\newcommand{\oper}[2]{\newcommand{#1}{\mathop{\mathrm{#2}}\nolimits} }
\oper{\DGA}{DGA}

\newcommand{\dbar}{\bar\partial}

\newcommand{\Lie}[1]{\operatorname{\textsl{#1}}}
\newcommand{\lie}[1]{{\operatorname{\mathfrak{#1}}}}

\newcommand{\End}{\Lie{End}}

\newcommand{\CC}{\mathbb{C}}

\newtheorem{thm}{Theorem}[section]
\newtheorem{lem}[thm]{Lemma}
\newtheorem{pro}[thm]{Proposition}
\newtheorem{cor}[thm]{Corollary}
\newtheorem{definition}[thm]{Definition}

\def\qed{\rule{2.3mm}{2.3mm}}
\newcommand{\bproof}{\noindent{\it Proof: }}
\newcommand{\eproof}{\hfill \qed \vspace{0.2in}}

\journal{a journal}

\begin{document}

\begin{frontmatter}



\title{Weak Mirror Symmetry \\ of Complex Symplectic
Algebras}


\author{R. Cleyton}
\ead{cleyton@gmail.com}
\address{Institute f\"ur Mathematik,
Humboldt Universit\"at zu Berlin, Unter den Linden 6, \\ D-10099 Berlin,
Germany}

\author{ Y. S. Poon} 
\ead{ypoon@ucr.edu} 
\address{ Department of Mathematics, University of California, Riverside, CA 92521,
USA}

\author{ G. P. Ovando \corref{cor1} \fnref{pa}}
\ead{ gabriela@fceia.unr.edu.ar} 
\cortext[cor1]{Corresponding author}
\fntext[pa]{Permanent address: FCEIA,
Universidad Nacional de Rosario, Av. Pellegrini 250, (2000) Rosario, Argentina.}
\address{Mathematisches Institut, Albert-Ludwigs-Universit\"at
Freiburg, Eckerstr. 1, \\ D-79104 Freiburg, Germany.}

\begin{abstract}
A complex symplectic structure on a Lie algebra $\lie h$ is an
integrable complex structure $J$ with a closed non-degenerate
$(2,0)$-form. It is determined by $J$ and the real part $\Omega$ of
the $(2,0)$-form. Suppose that $\lie h$ is a semi-direct product
$\lie g\ltimes V$, and both $\lie g$ and $V$ are Lagrangian with
respect to $\Omega$ and totally real with respect to $J$. This note
shows that $\lie g\ltimes V$ is its own weak mirror image in the
sense that the associated differential Gerstenhaber algebras
controlling the extended deformations of $\Omega$ and $J$ are
isomorphic.

The geometry of $(\Omega, J)$ on the semi-direct product $\lie
g\ltimes V$ is also shown to be equivalent to that of a torsion-free
flat symplectic connection on the Lie algebra $\lie g$.  By further
exploring a relation between $(J, \Omega)$ with  hypersymplectic
algebras, we find an inductive process to build families of complex
symplectic algebras of dimension $8n$ from the data of the
$4n$-dimensional ones.
\end{abstract}


\end{frontmatter}



\section{Introduction}

In this paper, we continue our analysis on weak mirror symmetry of
Lie algebras with a focus on complex symplectic structures.

Recall that the extended deformation theory of complex or symplectic
structures are dictated by  their associated differential
Gerstenhaber algebras (a.k.a. DGA) and the induced Gerstenhaber
algebra structure on its cohomology ring \cite{BK}  \cite{Zhou}. If
$M$ is a manifold with a complex structure $J$ and $M^\vee$ is
another manifold with a symplectic structure $\omega$ with the same
dimension, then $(M, J)$ and $(M^\vee, \omega)$ form a weak mirror
pair if the associated differential Gerstenhaber algebras $\DGA(M,
J)$ and $\DGA(M^\vee, \omega)$ are quasi-isomorphic \cite{Mer-note}.

Inspired by the concept of T-duality in mirror symmetry \cite{SYZ}
\cite{Ben}, we develop a set of tools to analyze weak mirror
symmetry on semi-direct product Lie algebras \cite{CLP}, and apply
it to study complex and symplectic structures on real
six-dimensional nilpotent algebras in details \cite{CP} \cite{CLP}.

However, the work in \cite{CLP} does not explain an elementary, but
non-trivial example on Kodaira-Thurston surface, a real
four-dimensional example \cite{Poon}. Yet Kodaira-Thurston surface
is a key example of hyper-symplectic structures \cite{Kamada}
\cite{Fino}. As explained by Hitchin, from a $G$-structure
perspective, hyper-symplectic structures are the non-compact
counterpart of hyper-K\"ahler structures. The unifying feature is a
complex symplectic structure, i.e. a $\dbar$-closed holomorphic
non-degenerate $(2,0)$-form on complex even-dimensional manifolds
\cite{Hyper-sym}.

In this notes, we analyze complex symplectic structures adopted to a
semi-direct product Lie algebra.

Our first observation is Proposition \ref{basic omega}, which
essentially states that if $\mathfrak g$ is a Lie algebra symplectic
structure $\omega$ and a torsion-free, flat, symplectic connection
$\gamma$, then the semi-direct product $\mathfrak g\ltimes_\gamma V$
admits a complex symplectic structure $\Omega$ such that $\mathfrak
g$ and $V$ are totally real, and Lagrangian with respect to the real
part of $\Omega$. Moreover, every complex symplectic structure on a
Lie algebra with such characteristics arises in exactly this way.
Algebras with such complex symplectic structures are called \it
special Lagrangian. \rm

With a computation using dual representation and the results in
\cite{CP} and \cite{CLP}, we derive the main result of this notes in
Theorem \ref{main}, which states that a semi-direct product Lie
algebra, special Lagrangian with respect to the complex symplectic
structure is its own weak mirror image.

It is known that hypersymplectic structures form a special sub-class
of complex symplectic structure. In Section \ref{sec: hyper}, we
explore this relation, and find in Theorem \ref{tower} that
whenever one gets a special Lagrangian complex symplectic Lie
algebra, one gets a tower of such objects in higher dimension.

 In the last section, we produce a few examples
explicitly. In particular, we put the work of the second author on
Kodaira-Thurston within the framework of complex symplectic
structures and provide a theoretical framework to explain the ad hoc
computation in \cite{Poon}. A set of new examples generalizing the
hypersymplectic structure on Kodaira-Thurston surfaces is
constructed.
\section{Summary of Weak Mirror Symmetry}
We briefly review the basic definitions and  results in \cite{CP}
and \cite{CLP}. Except for key concepts, details are referred to
these two papers.

It is known that for a manifold $M$ with a complex structure $J$
there associates a differential Gerstenhaber algebra $\DGA(M, J)$
over $\mathbb C$. This object dictates the deformation theory of the
complex structure. For a manifold $N$ with a symplectic form
$\Omega$, there also associates a differential Gerstenhaber algebra
structure $\DGA(N, \Omega)$. It dictates its deformation theory in
an extended sense \cite{BK}. If $M$ and $N$ are manifolds of the
same dimension, and if $\DGA(M, J)$ and $\DGA(N, \omega)$ are
quasi-isomorphic, then the pair of structures $(M, J)$ and $(N,
\Omega)$ are said to be a weak mirror pair \cite{CLP}
\cite{Mer-note}. We focus on the case when the manifolds in
questions are Lie groups, and all geometric structures are
left-invariant. Therefore, all constructions are expressed in terms
of Lie algebras. In particular, exterior differential becomes the
Chevalley-Eilenberg differentials on the dual of a Lie algebra. In
this notes, if  $\DGA(1)$ and $\DGA(2)$ are two quasi-isomorphic
differential Gerstenhaber algebras, we denote such a relation by
$\DGA(1)\sim\DGA(2)$. The notation $\DGA(1)\cong\DGA(2)$ denotes
isomorphism.

Suppose $\mathfrak h$ is the semi-direct product of a subalgebra $\mathfrak g$
and an abelian ideal $V$: $\mathfrak h=\mathfrak g\ltimes V$. The semi-direct
product structure is said to be totally real with respect to an
integrable complex structure $J$ on $\mathfrak h$ if
\[ J\mathfrak g=V \quad
\mbox{ and } \quad JV=\mathfrak g.
\]
  It is said to be Lagrangian with respect
to a symplectic form $\Omega$ on $\mathfrak h$ if both $\mathfrak g$ and $V$
are Lagrangian. By contraction in the first variable, a two-form on
$\mathfrak h$ is a linear map from $\mathfrak h$ to ${\mathfrak h}^*$. It is
non-degenerate if the induced map is a linear isomorphism. If the
semi-direct product is Lagrangian, then
\[
\Omega: \mathfrak g\to V^* \quad \mbox{ and }
\quad \Omega: V\to {\mathfrak g}^*.
\]

Given a semi-direct product $\mathfrak h=\mathfrak g\ltimes V$, the adjoint
action of $\mathfrak g$ on $V$ is a representation of $\mathfrak g$ on $V$.
Conversely, given any presentation $\gamma:\mathfrak g\to \End(V)$, one
could define a Lie bracket on the vector space $\mathfrak g\oplus V$ by
\begin{equation}\label{semi d}
[(x,u),(y,v)]=([x,y], \gamma(x)v-\gamma(y)u),
\end{equation}
where $x, y\in \mathfrak g$ and $u,v\in V$. The resulting Lie algebra is
denoted by $\mathfrak g\ltimes_\gamma V$. One could take the dual
representation $\gamma^*:\mathfrak g\to \End(V^*)$: for any $v$ in $V$,
$\alpha$ in $V^*$ and $x$ in $\mathfrak g$,
\[
(\gamma^*(x)\alpha)v:=-\alpha(\gamma(x)v).
\]
We address the Lie algebra ${\widehat{\mathfrak h}}:=\mathfrak
g\ltimes_{\gamma^*}V^*$ as the dual semi-direct product of $\mathfrak
h=\mathfrak g\ltimes_{\gamma}V$.

In \cite{CLP}, it is noted that weak mirror pairs on semi-direct
products are given in terms of dual semi-direct products.
\begin{pro}{\rm\cite{CLP}}\label{prop: cx to sym}
Let $\mathfrak h=\mathfrak g\ltimes_{\gamma}V$ be a semi-direct product
totally real with respect to an integrable complex structure $J$. On
the dual semi-direct product ${\widehat{\mathfrak h}}=\mathfrak
g\ltimes_{\gamma^*}V^*$, define
\begin{equation}\label{cx to sym}
{\widehat\Omega}((x, \mu), (y, \nu)):=\nu(Jx)-\mu(Jy),
\end{equation}
where $x,y\in \mathfrak g$, $\mu, \nu\in V^*$. Then $\widehat\Omega$ is a
symplectic form, and the dual semi-direct product $\widehat{\mathfrak h}$
is Lagrangian with respect to $\widehat\Omega$. Moreover, there is a
natural isomorphism $\DGA(\mathfrak h, J)\cong \DGA({\widehat{\mathfrak h}},
{\widehat\Omega})$.
\end{pro}

\begin{pro}{\rm\cite{CLP}}\label{prop: sym to cx}
Let $\mathfrak h=\mathfrak g\ltimes_{\gamma}V$ be a semi-direct product
Lagrangian with respect to a symplectic form $\Omega$. On the dual
semi-direct product ${\widehat{\mathfrak h}}=\mathfrak
g\ltimes_{\gamma^*}V^*$, define
\begin{equation}\label{sym to cx}
{\widehat J}(x,\mu):=(-\Omega^{-1}(\mu), \Omega(x)),
\end{equation}
where $x\in \mathfrak g$ and $\mu\in V^*$. Then $\widehat J$ is an
integrable complex structure, and the dual semi-direct product
$\widehat{\mathfrak h}$ is totally real with respect to $\widehat J$.
Moreover, there is a natural isomorphism $\DGA(\mathfrak h, \Omega)\cong
\DGA({\widehat{\mathfrak h}}, {\widehat J})$.
\end{pro}

\section{Complex Symplectic Algebras}

In this section, we present a few elementary facts on complex
symplectic structures on Lie algebras. One may globalize the
geometry to Lie groups by left translation. In favorable situation,
the geometry descends to a co-compact quotient. Therefore, lots of
consideration in this section has its counterpart on manifolds. A
discussion on the relation between complex symplectic structures,
hyperK\"ahler structures and hypersymplectic structures could be
found in \cite{Hyper-sym}. We stay focus on Lie algebras.

\begin{definition}{\rm \cite{Besse}, \cite{Hyper-sym}}
A complex symplectic structure on a real Lie algebra $\mathfrak h$ is a
pair $(J, \Omega_c)$ where $J$ is an integrable complex structure
and $\Omega_c$ is a non-degenerate closed  $(2,0)$-form.
\end{definition}

Let $\Omega_1$ and $\Omega_2$ be the real and imaginary parts of
$\Omega_c$ respectively, $\Omega_c=\Omega_1+i\Omega_2$. It is
immediate that they are closed real two-forms on $\mathfrak h$. Define
\[
\mathfrak h^{1,0}=\{X-iJX\in {\mathfrak h}_\CC: X\in \mathfrak h\}, \quad \mathfrak
h^{0,1}=\{X+iJX\in {\mathfrak h}_\CC: X\in \mathfrak h\}.
\]
They are respectively the $+i$ and $-i$ eigenspace of the complex
structure on the complexification of $\mathfrak h$.

Since the complex structure $J$ is integrable, $\dbar\Omega_c=0$ and
$\partial\Omega_c=0$. By definition of non-degeneracy, the map \[
\Omega_c: {\mathfrak h}^{1,0}\to {\mathfrak h}^{*(1,0)}
\]
is a linear isomorphism.

\begin{lem}\label{3.2} On a real vector space with an almost complex structure
$J$ the real part $\Omega_1$ and imaginary part $\Omega_2$ of a
non-degenerate complex $(2,0)$-form $\Omega_c$ are non-degenerate
two-forms with the properties
\begin{eqnarray}
&\Omega_2(X,Y)=-\Omega_1(JX, Y);&\label{Omega2}\\
&\Omega_1(JX,Y)=\Omega_1(X,JY), \quad \Omega_2(JX, Y)=\Omega_2(X,JY).
& \label{Omega1}
\end{eqnarray}

Conversely, given any real non-degenerate two-form $\Omega$
satisfying {\rm (\ref{Omega1})}, let $\Omega_1=\Omega$ and define
$\Omega_2$ by {\rm (\ref{Omega2})}, then
$\Omega_c=\Omega_1+i\Omega_2$  is a non-degenerate $(2,0)$-form.
\end{lem}
\bproof
Suppose $\Omega_c$ is a $(2,0)$-form on a real vector space $W$ with
almost complex structure $J$ for any real vectors $X$ and $Y$,
$\Omega_c(X+iJX, Y)=0.$ Setting $\Omega_c=\Omega_1+i\Omega_2$, and
taking the real and imaginary part, one get the relation in
(\ref{Omega2}).

Since $\Omega_1$ is real, $\Omega_1(X+iJX, Y+iJY)=\Omega_1(X-iJX,
Y-iJY)$. The  identity (\ref{Omega1}) for $\Omega_1$ follows. An
identical computation shows that $\Omega_2$ satisfies
(\ref{Omega1}).

Conversely, given a real two-form $\Omega_1$ satisfying
(\ref{Omega1}). Define $\Omega_2$ by (\ref{Omega2}), then $\Omega_2$
satisfies (\ref{Omega1}). It follows that
$\Omega_c:=\Omega_1+i\Omega_2$ is a $(2,0)$-form.

Consider the non-degeneracy of $\Omega_c$.  $\Omega_c(X-iJX)=0$ if
and only if
 \[
 \Omega_1(X)+\Omega_2(JX)=0 \quad \mbox{ and } \quad
 \Omega_1(JX)-\Omega_2(X)=0.
 \]
By (\ref{Omega2}), $\Omega_1(X)=\Omega_2(JX)$. Therefore,
$\Omega_c(X-iJX)=0$ if and only if $\Omega_1(X)=0$ and
$\Omega_2(X)=0$. Therefore, $\Omega_c$ is non-degenerate if and only
if both $\Omega_1$ and $\Omega_2$ are non-degenerate.
\eproof

On a real Lie algebra $\mathfrak h$ with complex structure $J$, a
$(2,0)$-form $\Omega_c=\Omega_1+i\Omega_2$ is closed if and only if
$\Omega_1$ and $\Omega_2$ are closed.

Given the complex structure, the Chevalley-Eilenberg differential
$d$ has a type decomposition. Denote its $(1,0)$-part by $\partial$.
Its $(0,1)$-part by $\dbar$. The integrability of $J$ implies that
$d\Omega_c=0$ if and only if $\partial\Omega_c=0$ and
$\dbar\Omega_c=0$.

Since $d\Omega_1=0$, for any $X, Y, Z$ in $\mathfrak h$,
\[
\dbar\Omega_1(X+iJX, Y-iJY, Z-iJZ )=0.
\]
As a consequence of the vanishing Nijenhuis tensor and
(\ref{Omega1}), the real and imaginary part of the above constraints
respectively
\begin{eqnarray}
&\Omega_1([X,Y],Z)+\Omega_1(J[JX,Y],Z)&\nonumber\\
&\hspace{1in} +\Omega_1([Z,X],Y)+\Omega_1(J[Z,JX],Y)=0.&
\label{real d omega}\\
&\Omega_1([X,Y],JZ)+\Omega_1(J[JX,Y],JZ)&\nonumber\\
&\hspace{1in} +\Omega_1([JZ,X],Y)+\Omega_1(J[JZ,JX],Y)=0.& \label{im
d omega}
\end{eqnarray}
The latter is apparently equivalent to the former as the map $J$ is
a real linear isomorphism from $\mathfrak h$ to itself.

\begin{lem}\label{o2 sym} Let $\mathfrak h$ be a Lie algebra with an integrable
complex structure $J$. Let $\Omega_1$ be a symplectic form on $\mathfrak
h$ satisfying {\rm (\ref{Omega1})}. Define $\Omega_2=-\Omega_1\circ
J$, then $\Omega_2$ is a symplectic form and satisfies {\rm
(\ref{Omega1})}.
\end{lem}
\bproof Let $X, Y, Z$ be any vectors in $\mathfrak h$. By
integrability of $J$,
\begin{eqnarray*}
&&d\Omega_2(X, Y, Z)=\\
&=&-\Omega_1(J[X, Y], Z)-\Omega_1(J[Y, Z], X)-\Omega_1(J[Z, X], Y)\\
&=&\Omega_1([JX, Y], Z)+\Omega_1([X, JY], Z)-\Omega_1([JX,JY], JZ)\\
&&\quad +\Omega_1([JY, Z], X)+\Omega_1([Y, JZ], X)-\Omega_1([JY,JZ], JX)\\
&& \quad \quad +\Omega_1([JZ,X],Y)+\Omega_1([Z, JX],
Y)-\Omega_1([JZ,JX], JY).
\end{eqnarray*}
Since $d\Omega_1=0$, the above is equal to
\begin{eqnarray*}
&&\Omega_1([JX, Y], Z)+\Omega_1([Z, JX], Y) -\Omega_1([Z, X],
JY)-\Omega_1([X,Y], JZ)\\
&=&-\Omega_1(J[JX, Y], JZ)-\Omega_1([X,Y], JZ)+\Omega_1([Z, JX], Y)
-\Omega_1([Z, X], JY).
\end{eqnarray*}
By (\ref{im d omega}), it is equal to
\[
\Omega_1([JZ, X]+J[JZ, JX]+[Z, JX]-J[Z, X], Y),
\]
which is equal to zero by integrability of $J$. \eproof

\begin{pro}\label{character} Given a complex structure $J$ on a Lie algebra $\mathfrak h$,
a two-form $\Omega_c=\Omega_1+i\Omega_2$ is non-degenerate and
closed if and only if $\Omega_1$ is a symplectic form satisfying
{\rm (\ref{Omega1})} and $\Omega_2=-\Omega_1\circ J$.

Conversely, if $\Omega_1$ is a symplectic form on $\mathfrak h$
satisfying {\rm (\ref{Omega1})}, define $\Omega_2:=-\Omega_1\circ
J$, then $\Omega_1+i\Omega_2$ is a complex symplectic structure.
\end{pro}

\bproof The first part of this proposition is a paraphrase of Lemma
\ref{3.2}.

Conversely, given $\Omega_1=\Omega$. Define $\Omega_2$ by
(\ref{Omega2}). Since the complex structure is integrable, and
$\Omega_1$ satisfies {\rm (\ref{Omega1})}, $d\Omega_1=0$ implies
that $\Omega_1$ satisfies (\ref{real d omega}). As a consequent of
Lemma \ref{o2 sym} $\Omega_2=-\Omega_1\circ J$ is a symplectic form.
\eproof

\section{Special Lagrangian Structures}

As a consequence of Proposition \ref{character}, whenever $J$ is a
complex structure and $\Omega$ is a \it real \rm symplectic form on
a Lie algebra $\mathfrak h$ such that (\ref{Omega1}) holds, the pair $(J,
\Omega)$ is addressed as a complex symplectic structure. However,
one should note that if $\Omega_c$ is a complex symplectic form with
respect to $J$, so is $re^{-i\theta}\Omega_c$ for any $r>0$ and
$\theta\in[0,2\pi)$. While the real factor $r$ will only change the
symplectic forms $\Omega_1$ and $\Omega_2$ by real homothety, the
factor $e^{-i\theta}$ changes $\Omega_1$ to $\cos\theta\
\Omega_1+\sin\theta\ \Omega_2$. Therefore, the geometry related to a
choice of $\Omega_1$ or $\Omega$ is not invariant of the complex
homothety on $\Omega_c$. As we choose to work with a particular
choice of $\Omega_1$ within complex homothety, straightly speaking,
we work with \it polarized \rm complex symplectic structures.

We examine the characterization of this type of objects when the
underlying algebra is a semi-direct product. A reason for working
with semi-direct product is due to the fact that the underlying
Gerstenhaber algebra $\DGA(\mathfrak h, J)$ of a Lie algebra $\mathfrak h$
with a complex structure $J$ is the exterior algebra of a
semi-direct product Lie algebra. Therefore, when one searches for a
weak mirror image for $(\mathfrak h, J)$, it is natural to search among
semi-direct products. Details for this discussion could be found in
\cite{CLP}. In our current investigation, the geometry and algebra
together also force upon us semi-direct product structures.

\begin{lem}\label{ideal} Let $\mathfrak h$ denote a  Lie algebra equipped
with a complex symplectic structure $(J, \Omega)$. Assume $V\subset
\mathfrak h$ is an isotropic ideal. Then $V$ is abelian; and $JV$ is a
subalgebra. Moreover if  $V$ is totally real and Lagrangian then
$\mathfrak h= JV \ltimes V $.
\end{lem}
\bproof Let $x,y\in V$.
The closeness condition of $\Omega$ says for all $z$ in $V$,
\[0=\Omega([x,y], z)+\Omega([y,z], x)+\Omega([z,x], y)=\Omega([x,y],
z).
\]
Since $\Omega$ is non-degenerate, $[x,y]=0$ for all $x,y \in V$. It
proves that $V$ is abelian.

Using the fact that $V$ is abelian and  the  integrability condition
of $J$ one gets
$$[Jx,Jy]=J([Jx,y]+[x,Jy]) \qquad \forall x,y \in V$$
which proves that $JV$ is a subalgebra.

If $V$ is Lagrangian, then $JV$ is also Lagrangian. Moreover
whenever $V$ is totally real then $V \cap J V=\{0\}$ and so one gets
 $\mathfrak h= JV\ltimes V$.
\eproof

\begin{definition} Let $\mathfrak h$ be a real Lie algebra, $V$ an
abelian ideal in $\mathfrak h$ and $\mathfrak g$ a complementary subalgebra so
that $\mathfrak h=\mathfrak g\ltimes V$. We say that a complex symplectic
structure $(J, \Omega)$ on $\mathfrak h$ is special Lagrangian with
respect to the semi-direct product if  $\mathfrak g$ and $V$ are totally
real with respect to $J$ and isotropic with respect to $\Omega$.
\end{definition}

The rest of this section is to characterize the geometry of a
special Lagrangian structure $(J, \Omega)$ on $\lie g\ltimes V$ in
terms of a connection on $\lie g$.

\

 The semi-direct product is given by the restriction
of the adjoint representation of $\mathfrak h$ to $\mathfrak g$. It
becomes a representation $\rho$ of the Lie subalgebra $\mathfrak g$
on the vector space $V$. There follows  a torsion-free flat
connection $\gamma$ associated to the complex structure $J$
\cite{CLP}. For any $x, y$ in $\mathfrak g$,
\begin{equation}\label{gamma}
\gamma(x)y:=-J\rho(x)Jy.
\end{equation}
Let $\Omega_1=\Omega$ be a two-form on $\mathfrak h$ in Proposition
\ref{character}. It satisfies (\ref{im d omega}). Apply this formula
to $x, y, Jz$ when $x,y,z\in \mathfrak g$. Since $Jz, Jx$ are in $V$, and
$V$ is abelian, $[Jz,Jx]=0.$  By definition of $\gamma$ and the
vanishing of its torsion,  Identity (\ref{im d omega}) is reduced to
\begin{eqnarray*}
0&=&\Omega(\gamma(x)y-\gamma(y)x,
Jz)+\Omega(\gamma(y)x, Jz)+\Omega(\gamma(x)z, Jy)\\
&=&\Omega(\gamma(x)y, Jz)-\Omega(\gamma(x)z, Jy).
\end{eqnarray*}
As $\Omega_2=-\Omega\circ J$. The above is equivalent to
\begin{equation}\label{closed on g}
\Omega_2(\gamma(x)y, z)-\Omega_2(\gamma(x)z,y)=\Omega_2(\gamma(x)y,
z)+\Omega_2(y, \gamma(x)z)=0.
\end{equation}
Note that $\Omega_2$ is a symplectic form on $\mathfrak h$. If
$\Omega_2(x,y)=\Omega_1(Jx, y)=0$ for all $x\in \mathfrak g$, then $y=0$
because $J\mathfrak g=V$ and $\Omega_1$ is isotropic on $\mathfrak g$.
Therefore, the restriction of $\Omega_2$ on $\mathfrak g$ is
non-degenerate. So is the restriction of $\Omega_2$ on $V$. It
follows that $\Omega_2$ is restricted to symplectic forms on $\mathfrak
g$ and $V$. Due to Identity (\ref{Omega1}) and the assumption on $J$
being totally real, the symplectic structure on $V$ is completely
dictated by the restriction of $\Omega_2$ on $\mathfrak g$. Denote this
restriction on $\mathfrak g$ by $\omega$. Then (\ref{closed on g}) is
equivalent to the symplectic form $\omega$ on $\mathfrak g$ being
parallel with respect to the connection $\gamma$. i.e.
\begin{equation}\label{parallel
omega} \omega(\gamma(x)y, z)+\omega(y, \gamma(x)z)=0.
\end{equation}

In summary, $\gamma$ is a torsion-free flat symplectic connection on
$\lie g$ with respect to the symplectic form $\omega$.

Conversely, suppose $\mathfrak g$ is equipped with a symplectic form
$\omega$, parallel with respect to a flat torsion-free connection
$\gamma$ on $\mathfrak g$. Let $V$ be the underlying vector space of
$\mathfrak g$. Define $\mathfrak h:=\mathfrak g\ltimes_\gamma V$,
set
\begin{equation}\label{the j}
J(x,y)=(-y, x),
\end{equation}
 then $J$ is an integrable complex
structure.  $\lie g$ and $V$ are totally real. Note that in this
case, for $x,y\in \mathfrak g$,
\begin{eqnarray*}
&&-J\gamma(x,0)J(y,0) =-J\gamma(x,0)(0, y) =-J(0, \gamma(x)y)
=-(-\gamma(x)y,0)\\ &=&(\gamma(x)y,0).
\end{eqnarray*}
If we use the notation $x$ to denote $(x,0)$ in $\mathfrak h$, then the
above identity becomes
\[
-J\gamma(x)Jy=\gamma(x)y.
\]
Therefore, one may consider the representation $\gamma$ playing both
the role of representation $\rho$ and the role of the connection
$\gamma$ in (\ref{gamma}). This coincidence is the consequence of
$V$ being the underlying vector space of the Lie algebra $\mathfrak
g$.

Next, define
\begin{eqnarray}
&\Omega_1((x,u), (y,v)):=-\omega(x,v)-\omega(u,y).& \label{o1 formula}\\
&\quad\quad \Omega_2((x,u), (y,v)):=\omega(x,y)-\omega(u, v).&
\label{o2 formula}
\end{eqnarray}
It is apparent that $\Omega_1$ and $\Omega_2$ are non-degenerate
skew-symmetric form on $\mathfrak h$, with $\lie g$ and $V$
isotropic with respect to $\Omega_1$.  Moreover,
\begin{eqnarray*}
&&\Omega_1(J(x,y), J(x',y'))=\Omega_1((-y,x), (-y',x'))
=-\omega(-y,x')-\omega(x,-y')\\
&=&-\Omega_1((x,y), (x',y')).\\
&&\Omega_2(J(x,y), J(x',y'))=\Omega_2((-y,x),
(-y',x'))=\omega(-y,-y')-\omega(x,x')\\
&=&-\Omega_2((x,y), (x',y')).\\
&&-\Omega_1((x,y), J(x',y'))=-\Omega_1((x,y), (-y',x'))=\omega(x,
x')+\omega(y,-y')\\
&=&\Omega_2((x,y), (x',y')).
\end{eqnarray*}
Therefore, (\ref{Omega2}) and (\ref{Omega1}) are satisfied.
\begin{eqnarray*}
&&-d\Omega_1((x_1,y_1),(x_2,y_2), (x_3,y_3))\\
&=&\Omega_1([(x_1,y_1),(x_2,y_2)], (x_3,y_3))+\mbox{cyclic permutation}\\
&=&\Omega_1(([x_1, x_2],\gamma(x_1)y_2-\gamma(x_2)y_1),
(x_3,y_3))+\mbox{cyclic permutation}\\
&=&\omega([x_1, x_2], y_3)+\omega(\gamma(x_1)y_2-\gamma(x_2)y_1,
x_3)+\mbox{cyclic permutation}\\
&=&\omega(\gamma(x_1)x_2, y_3)-\omega(\gamma(x_2)x_1, y_3)+
\omega(\gamma(x_1)y_2, x_3)-\omega(\gamma(x_2)y_1, x_3)\\
&&+\mbox{cyclic permutation}.
\end{eqnarray*}
It is equal to zero because $\omega$ is parallel with respect to
$\gamma$ (\ref{parallel omega}). Therefore, $\Omega_1$ is
symplectic. $\Omega_2$ is also a symplectic form due to $d\omega=0$
and (\ref{parallel omega}). Therefore, $\Omega_c=\Omega_1+i\Omega_2$
is closed, and it defines a complex symplectic structure. Below is a
summary of our observations in this section.

\begin{pro}\label{basic omega} Let $\mathfrak h=\mathfrak g\ltimes_\gamma V$
be a semi-direct product algebra. Suppose that it admits
$(\Omega_1+i\Omega_2, J)$ as a special Lagrangian complex symplectic
structure. Let $\omega$ be the restriction of $\Omega_2$ on $\mathfrak
g$. Then $\omega$ is a symplectic structure on the Lie algebra $\mathfrak
g$, $\gamma$ is a torsion-free, flat connection on $\mathfrak g$
preserving the symplectic structure $\omega$.

Conversely, given a Lie algebra $\mathfrak g$ with a symplectic
structure $\omega$ and a torsion-free flat symplectic connection,
the semi-direct product $\mathfrak g\ltimes_\gamma V$ admits a
complex symplectic structure $(J, \Omega)$ such that both $\lie g$
and $V$ are totally real with respect to $J$ and Lagrangian with
respect to $\Omega$.
\end{pro}

\section{The Dual Semi-direct Product}

Let $\mathfrak h=\mathfrak g\ltimes_\gamma V$ be a semi-direct product with a
special Lagrangian complex holomorphic structure $(J, \Omega)$.
Recall that $\gamma:\mathfrak g\to \End (V)$ is a representation. It
induces a dual representation: $\gamma^*:\mathfrak g\to \End(V^*)$.

Use $x,y,z$ and $a, b, c$ to represent generic elements in
$\mathfrak g$. $\alpha, \beta$ to represent generic elements in
${\mathfrak g}^*$, $u, v, w$ to represent generic elements in $V$,
$\mu,\nu,\zeta$ to represent generic elements in $V^*$. By
Proposition \ref{basic omega}, the restriction of
\begin{equation}\label{little o}
\omega(x,y):=-\Omega(Jx,y)
\end{equation}
 is a symplectic form on
$\mathfrak g$. By contractions, $\omega$ defines an isomorphism.

\begin{lem}\label{lem:iso} Let $V$ be the underlying symplectic vector space of Lie algebra
$\mathfrak g$. The linear map $\varpi: \mathfrak g\ltimes_\gamma V\to \mathfrak g
\ltimes_{\gamma^*} V^*$ defined by
\begin{equation}
\varpi(x,u):=(x, \omega(u))
\end{equation}
is a Lie algebra isomorphism.
\end{lem}
\bproof By definition of the semi-direct product $\mathfrak
g\ltimes_{\gamma^*}V^*$,
\[
[(x,\omega(u)), (y, \omega(v))]=([x,y],
\gamma^*(x)\omega(v)-\gamma^*(y)\omega(u)).
\]
For all $c\in V$,
\begin{eqnarray*}
&&\Big(\gamma^*(x)\omega(v)-\gamma^*(y)\omega(u)\Big)(c)\\
&=&-\omega(v)(\gamma(x)c)+\omega(u)(\gamma(y)c)=-\omega(v,
\gamma(x)c)+\omega(u, \gamma(y)c).
\end{eqnarray*}
Since $\gamma$ is a symplectic connection,  the last quantity is
equal to
\begin{eqnarray*}
&&\omega(\gamma(x)v, c)-\omega(\gamma(y)u, c)
=\omega(\gamma(x)v-\gamma(y)u, c)\\
&=&\omega(\gamma(x)v-\gamma(y)u)c.
\end{eqnarray*}
As $[(x,u), (y,v)]=([x,y], \gamma(x)v-\gamma(y)u)$,
\[
[(x, \omega(u)), (y, \omega(v))]=([x, y],
\omega(\gamma(x)v-\gamma(y)u)).
\]
\eproof

On $\mathfrak g\ltimes_{\gamma^*}V^*$, define $\widetilde J$ and
$\widetilde\Omega$ by
\begin{equation}\label{tilded}
{\widetilde J}:=\varpi\circ J\circ \varpi^{-1}, \quad
\varpi^*{\widetilde \Omega}=\Omega.
\end{equation}
Since $\varpi$ is an isomorphism of Lie algebra, $\widetilde J$ and
$\widetilde\Omega$ defines a complex symplectic structure on
$\widehat{\mathfrak h}=\mathfrak g\ltimes_{\gamma^*}V^*$. Explicitly,
\begin{eqnarray}
&{\widetilde J}(x,\omega(u))=\varpi(J(x,u))=\varpi(Ju,Jx)
=(Ju, \omega(Jx)).&\label{tilded j}\\
 &{\widetilde\Omega}((x,\omega(u)),
(y,\omega(v)))=\Omega((x,u), (y,v))=\Omega(x,v)+\Omega(u,
y).&\label{tilded o}
\end{eqnarray}
Also, define ${\widetilde\Omega}_2:=-{\widetilde\Omega}\circ
{\widetilde J}$. By (\ref{tilded j}) and (\ref{tilded o}),
\begin{eqnarray*}
&&{\widetilde\Omega}_2((x,\omega(u)), (y,\omega(v)))
=-{\widetilde\Omega}_1({\widetilde J}(x,\omega(u)), (y,\omega(v)))\\
&=&-{\widetilde\Omega}((Ju,\omega(Jx)), (y,\omega(v)))=-\Omega(Ju, v)-\Omega(Jx, y)\\
&=&-\Omega(J(x,u), (y,v))=\Omega_2((x,u), (y,v)).
\end{eqnarray*}

The restriction of ${\widetilde\Omega}_2$ on the  summand  $\mathfrak g$
in $\mathfrak g\ltimes_{\gamma^*} V^*$ is identical to the restriction of
the $\Omega_2$ on $\mathfrak g\ltimes_{\gamma} V$. It is the symplectic
two-form $\omega$.

By Proposition \ref{prop: cx to sym}, given a totally real complex
structure $J$ on $\mathfrak h=\mathfrak g\ltimes_\gamma V$, there exists a
symplectic form $\widehat\Omega$ on $\widehat{\mathfrak h}$. It is
determined by (\ref{cx to sym}). By (\ref{little o}),
$\omega(x)=-\Omega(Jx)$. i.e. $\omega(Jx)=\Omega(x)$.
\begin{eqnarray*}
&&{\widehat\Omega}((x,\omega(u)), (y,\omega(v)))=\omega(v)(Jx)-\omega(u)(Jy)\\
&=&\omega(v, Jx)-\omega(u, Jy)=-\Omega(Jv, Jx)+\Omega(Ju, Jy)\\
&=&\Omega(v,x)-\Omega(u,y).
\end{eqnarray*}
Comparing with (\ref{tilded o}), we find that
$\widehat\Omega=-\widetilde\Omega$. Now  we obtain isomorphisms of
differential Gerstenhaber algebras:
\begin{eqnarray*}
\DGA(\mathfrak h, J)&\cong& \DGA(\widehat{\mathfrak h}, \widehat{\Omega})
\quad \mbox{ by Proposition \ref{prop: cx to sym} }\\
&=&\DGA(\widehat{\mathfrak h}, \widetilde{\Omega}) \quad \mbox{ due to }
\quad \widehat\Omega=-\widetilde\Omega, \\
&\cong& \DGA({\mathfrak h}, \Omega) \quad \mbox{ by (\ref{tilded}) and
Lemma \ref{lem:iso}}.
\end{eqnarray*}

Similarly, by Proposition \ref{prop: sym to cx}, if the semi-direct
product $\mathfrak h$ is Lagrangian with respect to the symplectic form
$\Omega$, there exists a complex structure $\widehat J$ on
$\widehat{\mathfrak h}$. It is determined by (\ref{sym to cx}).
\begin{eqnarray*}
&&{\widehat J}(x,\omega(u))=(-\Omega^{-1}(\omega(u)), \Omega(x))
=(\Omega^{-1}(\Omega (Ju)), \omega(Jx))\\
 &=&(Ju,\omega(Jx)).
\end{eqnarray*}
Comparing with (\ref{tilded j}), we find that ${\widehat
J}={\widetilde J}$. Therefore, we again obtain the isomorphisms:
\begin{eqnarray*}
\DGA(\mathfrak h, \Omega)&\cong& \DGA(\widehat{\mathfrak h}, {\widehat J})
\quad \mbox{ by Proposition \ref{prop: sym to
cx}}\\
&=&\DGA(\widehat{\mathfrak h}, {\widetilde J}) \quad \mbox{ due to }
\quad {\widehat J}={\widetilde J}\\
&\cong& \DGA(\mathfrak h, J) \quad \mbox{ by (\ref{tilded}) and Lemma
\ref{lem:iso}}.
\end{eqnarray*}
The last two strings of isomorphisms yield the main observation in
this notes.

\begin{thm}\label{main} Every special Lagrangian complex symplectic
semi-direct product is isomorphic to its mirror image.

To be precise, let $\mathfrak h=\mathfrak g\ltimes V$ be a semi-direct
product, $J$ an integrable complex structure,
$\Omega_c=\Omega_1+i\Omega_2$ a closed $(2,0)$-form such that $\mathfrak
h$ is totally real with respect to $J$ and Lagrangian with respect
to $\Omega=\Omega_1$, then $\mathfrak h\cong {\widehat{\mathfrak h}}$ as Lie
algebras, and $ \DGA(\mathfrak h, J)\cong \DGA(\mathfrak h, \Omega)\cong \DGA(
{\widehat{\mathfrak h}}, {\widehat J})\cong \DGA({\widehat{\mathfrak h}},
{\widehat\Omega}). $
\end{thm}

\section{Relation to Hypersymplectic Geometry}\label{sec: hyper}

Given a special Lagrangian complex symplectic Lie algebra $\mathfrak
h=\mathfrak g\ltimes_\rho V$ with $(\Omega, J)$,  there are more
structures naturally tied to the given ones. For instance, consider
the linear map $E: \mathfrak g\ltimes_\rho V\to \mathfrak g\ltimes_\rho V$
defined by $E\equiv {\mathbf 1}$ on $\mathfrak g$ and $E\equiv -{\mathbf
1}$ on $V$. Then $E$ is an example of an {\it almost product
structure} \cite{AS}. Since both $\mathfrak g$ and $V$ are Lie
subalgebras of $\mathfrak h$, this almost product structure is
integrable, or a product structure without qualification. Since
$JV=\mathfrak g$ and $J\mathfrak g=V$,  $JE=-EJ$. Therefore, the pair $\{J,
E\}$ forms an example of a complex product structure \cite{AS}.

Given any complex product structure, there exists a unique
torsion-free connection $\Gamma$ on $\mathfrak h$ with respect to which
$J$ and $E$ are parallel. For $x,y,x',y' $ in $\mathfrak g$, define
\begin{equation}
\Gamma_{(x+Jx')}(y+Jy'):=-J\rho(x)Jy+ \rho(x)Jy'.
\end{equation}
In terms of the direct sum decomposition $\mathfrak h=\mathfrak g\ltimes_\rho
V$, we use $(x,0)$ to denote an element $x$ in $\mathfrak g$. Then for
any $y\in \mathfrak g$, $Jy=(0,y)$. It follows that the connection with
respect to the direct sum decomposition is given by
\begin{equation}\label{nabla}
\Gamma_{(x, x')}(y, y'):=(\rho(x)y,  \rho(x)y').
\end{equation}

 We check that this is a
torsion-free connection.
\begin{eqnarray*}
&&\Gamma_{(x+Jx')}(y+Jy')-\Gamma_{(y+Jy')}(x+Jx')
-[x+Jx', y+Jy']\\
&=&(-J\rho(x)Jy+ \rho(x)Jy')-(-J\rho(y)Jx+ \rho(y)Jx') - ([x,y]+\rho(x)Jy'-\rho(y)Jx')\\
&=&-J\rho(x)Jy+J\rho(y)Jx-[x,y]\\
&=&\gamma(x)y-\gamma(y)x-[x,y].
\end{eqnarray*}
Therefore, $\Gamma$ is torsion-free if and only if the connection
$\gamma$ is torsion-free.

To compute $\Gamma J$, we find
\begin{eqnarray*}
&&\Gamma_{(x+Jx')}J(y+Jy')-J\Gamma_{(x+Jx')}(y+Jy')\\
&=&\Gamma_{(x+Jx')}(-y'+Jy)-J(-J\rho(x)Jy+\rho(x)Jy')\\
&=&J\rho(x)Jy'+\rho(x)Jy-\rho(x)Jy-J\rho(x)Jy'=0.
\end{eqnarray*}
Therefore, $\Gamma J=0$. Similarly, $\Gamma E=0$ because
\begin{eqnarray*}
&&\Gamma_{(x+Jx')}E(y+Jy')-E\Gamma_{(x+Jx')}(y+Jy')\\
&=&\Gamma_{(x+Jx')}(y-Jy')-E(-J\rho(x)Jy+\rho(x)Jy')\\
&=&-J\rho(x)Jy-\rho(x)Jy'+J\rho(x)Jy-\rho(x)Jy'=0.
\end{eqnarray*}

Note that $\Gamma E=0$ implies that the connection $\Gamma$
preserves the eigenspaces $\mathfrak g$ and $V$. Therefore, one may
consider the restriction of $\Gamma$ onto the respective subspaces.
Indeed, the restriction of $\Gamma$ onto the $+1$ eigenspace
$\mathfrak g$ is precisely the connection $\gamma$ as seen in
(\ref{gamma}).

Recall that the real and imaginary parts of the complex symplectic
form are two different real symplectic forms. In a direct sum
decomposition, they are respectively given by (\ref{o1 formula}) and
(\ref{o2 formula}). We now consider a third differential form.
\begin{equation}\label{o3 formula}
\Omega_3((x,u),(y,v))=\omega(x,y)+\omega(u,v).
\end{equation}
The three differential forms are related by the endomorphisms $J$
and $E$:
\begin{eqnarray}
\Omega_1((x,u),(y,v))&=&\Omega_2(J(x,u), (y,v)), \\
\Omega_3((x,u),(y,v))&=&\Omega_2(E(x,u), (y,v)).
\end{eqnarray}
Treating the contraction with the first entry in any 2-form as an
endomorphism from a vector space to its dual, we have
\begin{eqnarray}
&\Omega_1=\Omega_2\circ J, \quad \Omega_3=\Omega_2\circ E, \quad
\Omega_1=\Omega_3\circ J\circ E,&\\
&\Omega_2=-\Omega_1\circ J, \quad \Omega_2=\Omega_3\circ E, \quad
\Omega_3=-\Omega_1\circ E\circ J.&
\end{eqnarray}

Define a bilinear form $g$ on $\mathfrak h$ by
\begin{equation}
g((x,u),(y,v)):=\Omega_2(JE(x,u), (y,v))=-\omega(x,v)+\omega(u,y).
\end{equation}
It is a non-degenerate symmetric bilinear form on a $4n$-dimensional
vector space with signature $(2n,2n)$, also known as a neutral
metric. Since
\begin{equation}
J\circ J=-{\mathbf I}, \quad E\circ E={\mathbf I}, \quad (J\circ
E)\circ (J\circ E)={\mathbf I}, \quad J\circ E=-J\circ E,
\end{equation}
where  ${\mathbf I}$ is the identity map, we could express the three
2-forms in terms of the neutral metric $g$ and the complex product
structures $\{J, E, JE\}$.
\begin{equation}
\Omega_1(\bullet, \bullet)=g(E\bullet, \bullet), \quad
\Omega_2(\bullet, \bullet)=g(JE\bullet, \bullet), \quad
\Omega_3(\bullet, \bullet)=g(J\bullet, \bullet).
\end{equation}

By \cite[Proposition 5]{AD}, $d\Omega_2=0$ implies that
$d\Omega_3=0$. Also $d\Omega_2=0$ if and only if $d\Omega_1=0$.
Given the complex symplectic structure, $\Omega_1$ and $\Omega_2$
are closed. Therefore, $\Omega_3$ is also closed. It follows that
$(g, \Omega_1, \Omega_2, \Omega_3)$ forms a hypersymplectic
structure. In addition, by \cite[Corollary 10]{AD}, the connection
$\Gamma$ is the Levi-Civita connection for the neutral metric $g$.
Since $\Gamma J=0$ and $\Gamma E=0$, it follows that
\begin{equation}
\Gamma \Omega_1=\Gamma\Omega_2=\Gamma\Omega_3=0.
\end{equation}

In addition, since $\Gamma_{Jx}=0$ for all $x\in \mathfrak g$, and it is
only a matter of tracing definitions to verify that
\[
\Gamma_x\Gamma_{x'}-\Gamma_{x'}\Gamma_x-\Gamma_{[x,x']}=0.
\]
The connection $\Gamma$ is flat. As a consequence of the second part
of Proposition \ref{basic omega}, we have the following.

\begin{thm}\label{tower} Let $\mathfrak g$ be a symplectic Lie algebra with a
torsion-free, flat symplectic connection $\gamma$ on the underlying
vector space $V$ of the Lie algebra. Then the space $\mathfrak h=\mathfrak
g\ltimes_\gamma V$ admits a hypersymplectic structure $\{\Omega_1,
\Omega_2, \Omega_3\}$ such that the Levi-Civita connection $\Gamma$
of the associated neutral metric is flat, and symplectic with
respect to each of the three given symplectic structures.

In particular, each of these three symplectic structures induces a
complex symplectic structure on $\mathfrak h\ltimes_\Gamma W$ where
$W$ is the underlying vector space of the algebra $\mathfrak h$.

\end{thm}

\section{Examples}

In this section, we present explicit computations. Theorem
\ref{tower} of the previous section provide further examples for
exploration.

To recap, $(\mathfrak g, \omega)$ denote a real Lie algebra equipped
with a symplectic structure. Let $V$ denote the underlying vector
space of $\mathfrak g$. We seek a linear map $\gamma: \mathfrak g\to
\End(V)$ such that for all $x,y, z\in \mathfrak g$, it is
\begin{itemize}
\item torsion-free: $\gamma(x) y-\gamma(y) x= [x,y]$;
\item symplectic: $\omega(\gamma(x) y, z)+\omega(y, \gamma(x)
z)=0$;
\item flat: $\gamma([x,y])=\gamma (x) \gamma(y)-\gamma(y) \gamma(x)$.
\end{itemize}
The last two conditions are equivalent to require the map $\gamma$
to admit a symplectic representation of $\mathfrak g$ on $(V, \omega)$:
\begin{equation}
\gamma: \mathfrak g\to \mathfrak{sp}(V,\omega).
\end{equation}

On $\mathfrak g$, choose a basis $\{e_1, \dots, e_{2n}\}$. Then the dual
basis is denoted by $\{e^1, \dots, e^{2n}\}$ such that the
symplectic form is
\[
\omega=e^1\wedge e^2+\cdots +e^{2n-1}\wedge e^{2n}.
\]
Let $V$ be the underlying symplectic vector space for the Lie
algebra $\mathfrak g$. The corresponding basis and dual basis are
respectively $\{v_1, \dots, v_{2n}\}$ and $\{v^1, \dots, v^{2n}\}$.
Let $\gamma_{jk}^{\ l}$ and $c_{jk}^{\ l}$ be the structural
constants for the representation $\gamma$ and the algebra $\mathfrak g$:
\[
\gamma(e_j)e_k=\gamma_{jk}^{\ l}e_l, \quad \mbox{ and } \quad [e_j,
e_k]=c_{jk}^{\ l}e_l.
\]
Then the structure equations for the semi-direct product $\mathfrak
h=\mathfrak g\ltimes_\gamma V$ are
\begin{equation}\label{structure}
[e_j, e_k]=c_{jk}^{\ l}e_l, \quad \quad [e_j, v_k]=\gamma_{jk}^{\
 l}v_l.
\end{equation}
The complex structure $J$ is determined by $Je_j=v_j$ as dictated by
(\ref{the j}). The three real symplectic forms on the semi-direct
product $\mathfrak h$ are given by (\ref{o1 formula}), (\ref{o2 formula})
and (\ref{o3 formula}) respectively:
\begin{eqnarray}
&\Omega_1=\sum_{j=1}^n(-e^{2j-1}\wedge v^{2j}+e^{2j}\wedge
v^{2j-1});&\label{o1}\\
&\Omega_2=\sum_{j=1}^n(e^{2j-1}\wedge e^{2j}-v^{2j-1}\wedge
v^{2j}).&\label{o2}\\
&\Omega_3=\sum_{j=1}^n(e^{2j-1}\wedge e^{2j}+v^{2j-1}\wedge
v^{2j}).&\label{o3}
\end{eqnarray}

\subsection{Four-dimensional examples}\label{4 d} \

There are only two two-dimensional Lie algebras, the abelian one
${\mathbb R}^2$ and the algebra $\mathfrak{aff}(\mathbb R)$ of
affine transformations on $\mathbb R$.

 Consider $\mathbb R^2$ as an abelian
Lie algebra with its canonical symplectic structure
$\omega=e^1\wedge e^2$.

\begin{pro}\label{hei}
 {\rm \cite{An}} All torsion-free, flat, symplectic connection on the
abelian algebra ${\mathbb R}^2$ is equivalent to $({\mathbb R}^2,
\gamma, e^1\wedge e^2)$, where $\gamma$ is either trivial {\rm
($\gamma=0$)}, or \[  \gamma(e_1)=\left(
\begin{matrix}
0 & 0 \\ 1 & 0
\end{matrix} \right), \quad \gamma(e_2)=0.
\]
\end{pro}

\

As a consequence, the resulting four-dimensional complex symplectic
Lie algebras are either  abelian or isomorphic to the trivial
one-dimensional extension of the three-dimensional Heisenberg Lie
algebra. We make further change of notations to facilitate future
computations. With respect to a vector space decomposition of $\mathfrak
g={\mathbb R}^2\ltimes_\gamma {\mathbb R}^2$, a basis is given by
\[
f_1=(e_1, 0)=e_1, \quad f_2=(e_2, 0)=e_2, \quad f_3=(0, e_1)=v_1,
\quad f_4=(0, e_2)=v_2.
\]
Then the non-trivial example given in Proposition \ref{hei} above
implies that $\gamma(e_1)f_3=f_4$ on the semi-direct product. The
structure equation for $\mathfrak g$ is given by
\begin{equation}
[f_1, f_3]=(0, \gamma(e_1)e_1)=(0, e_2)=f_4.
\end{equation}
The two symplectic forms on $\mathfrak g$ are respectively
\begin{equation}\label{ex: o1}
\Omega_1=-f^1\wedge f^4+f^2\wedge f^3, \quad \Omega_2=f^1\wedge
f^2-f^3\wedge f^4.
\end{equation}
Moreover, the complex structure $J$ on $\mathfrak g$ is
\begin{equation}
Jf_1=f_3, \quad Jf_2=f_4.
\end{equation}
Here we recover the complex and symplectic structures studied in
\cite{Poon}, and hence the related claim on the self-mirror property
\cite[Theorem 19]{Poon}.

In the context of hypersymplectic structure, the third 2-form is
\begin{equation}\label{ex: o3}
\Omega_3=f^1\wedge f^2+f^3\wedge f^4.
\end{equation}

 \

Consider next $\mathfrak{aff}(\mathbb R)$ the Lie algebra spanned by
$e_1,e_2$ with the Lie bracket $[e_1, e_2]=e_2$, with symplectic
structure $\omega=e^1\wedge e^2$.

\begin{pro}\label{solvable} {\rm \cite{An}}
All torsion-free, flat, symplectic connection on the algebra of
group of affine transformations  $\mathfrak{aff}(\mathbb R)$ is
equivalent to $(\mathfrak{aff}(\mathbb R), \gamma, e^1\wedge e^2)$,
where $\gamma$ is one of the following.
\begin{itemize}
\item $\gamma(e_1)= \left(
\begin{matrix}
-1 & 0 \\ 0 & 1
\end{matrix} \right)$, \quad   $\gamma(e_2)=0$.
\item
$\gamma(e_1)=\left(
\begin{matrix}
-1/2 & 0 \\ 0 & 1/2
\end{matrix} \right)$,  $\gamma(e_2)=\left(
\begin{matrix}
0 & 0 \\ -1/2 & 0
\end{matrix} \right)
$.
\end{itemize}
\end{pro}
With these connections on $\mathfrak{aff}(\mathbb R)$, one applies
Proposition \ref{basic omega} to construct two complex symplectic
algebras in real-dimension four.

\begin{cor}\label{cor: 4d} Suppose that $\mathfrak g\ltimes_\gamma V$ is a
non-abelian four-dimensional semi-direct product such that it is
special Lagrangian with respect to $(J, \Omega)$, a complex
symplectic structure, then it is one of three cases. With respect to
a basis $\{f_1, f_2, f_3,f_4\}$, the structure equations are
respectively
\begin{itemize}
\item $[f_1, f_3]=f_4$,
\item $[f_1, f_2]=f_2$, $[f_1, f_3]=-f_3$, $[f_1, f_4]=f_4$.
\item $[f_1, f_2]=f_2$, $[f_1, f_3]=-\frac12f_3$, $[f_1,
f_4]=\frac12 f_4,$ $[f_2, f_3]=-\frac12 f_4$.
\end{itemize}
The complex structure is determined by $Je_1=e_3$, $Je_2=e_4$. The
complex symplectic form is
$
\Omega_c=\Omega_1+i\Omega_2=i(f^1+if^3)\wedge (f^2+if^4).
$
where $ \Omega=\Omega_1=-f^1\wedge f^4+f^2\wedge f^3$ and
$\Omega_2=f^1\wedge f^2-f^3\wedge f^4$.
\end{cor}

\subsection{An example of Theorem \ref{tower}}
To construct an explicit example for Theorem \ref{tower}, we
consider the canonically induced structure on the non-trivial
semi-direct product extending $\mathfrak h:={\mathbb R}^2\ltimes_\gamma
{\mathbb R}^2$ as given in Proposition \ref{hei}. On this space, we
choose a basis in terms of the direct sum $\mathfrak h\oplus W$ where $W$
is the underlying vector space of $\mathfrak h$:
\begin{equation}
f_j=(f_j, 0), \quad v_{j}=(0, f_j), \quad \mbox{ for } 1\leq j\leq
4.
\end{equation}
By (\ref{nabla}), a torsion-free, flat, symmetric connection
$\Gamma$ on the direct sum $\mathfrak h\oplus W$ is given by
\begin{equation}
\Gamma_{f_j}=\left( \begin{matrix} \gamma(f_j) & 0 \\ 0 &
\gamma(f_j)
\end{matrix}\right), \quad
\Gamma_{v_{j}}=0.
\end{equation}
Since the only non-trivial term is given by $\gamma(f_1)f_3=f_4$,
the non-trivial terms for $\Gamma$ are given by
\begin{equation}
\Gamma_{f_1}f_3=f_4, \quad \Gamma_{f_1}v_3=v_4.
\end{equation}
They also determine the structural equations for the algebra $\mathfrak
h\ltimes_\Gamma W$:
\begin{equation}
[f_1, f_3]=f_4, \quad [f_1, v_3]=v_4.
\end{equation}
This algebra is eight-dimensional. Its center is
five-dimensional, spanned by $f_2, f_4, v_1, v_2, v_4$.
Since the commutator $C(\mathfrak h)=span\{f_4, v_4\}$ is
contained in the center, $\mathfrak h$ is a 2-step nilpotent Lie algebra.

Work in Section \ref{sec: hyper} shows that the connection $\Gamma$
is torsion-free, flat and symplectic with respect to the symplectic
forms $\Omega_1$, $\Omega_2$ in (\ref{ex: o1}) and $\Omega_3$ in
(\ref{ex: o3}). By Theorem \ref{tower}, any real linear combination
of these symplectic form determines a complex symplectic structure
on the real eight-dimensional algebra $\mathfrak h\ltimes_\Gamma W$.

\subsection{2-step nilmanifolds with hypersymplectic structures}
We shall produce 2-step nilpotent Lie algebras $\mathfrak h$
carrying hypersymplectic structures, generalizing the first Lie
algebra in Corollary \ref{cor: 4d}.

Let $\omega$ denote the canonical symplectic structure on $\mathbb R^{2n}$:
\begin{equation}\label{om}
\omega=\sum_{i=1}^n e^i\wedge e^{i+n}.
\end{equation}
Let $A_1, \dots, A_n$ be $n\times n$-symmetric matrices. Let
$\gamma:\mathbb R^{2n} \to \mathfrak{sp}(2n, \omega)$ be the linear
map  represented, in the  basis $e_1, e_2, \hdots, e_{2n}$, by the
matrices:
\[
\gamma(e_i)=\left(
\begin{matrix}
0 & 0 \\
A_i & 0
\end{matrix}
\right), \qquad  \gamma(e_{i+n})=0, \qquad  \mbox{ for all } i=1, 2,
\hdots n.
\] Treating $\mathbb R^{2n}$ as the trivial abelian
algebra, then $\gamma$ is a representation. Notice that
$[\gamma(e_i), \gamma(e_j)]=0$ so that $\gamma$  gives rise to a
flat, symplectic connection on $\mathbb R^{2n}$.

The torsion-free condition on $\mathbb R^{2n}$ requires $\gamma(e_i)
e_j=\gamma(e_j)e_i$ for all $ i,j=1, 2, \hdots, n$. It is equivalent
to require that the $j$-th column of $A_i$ is equal to the $i$-th
column of $A_j$.  Hence for $n=1$ one gets
\[
\gamma(e_1)=\left(
\begin{matrix}
0 & 0 \\
a & 0 \\
\end{matrix}
\right), \qquad \gamma(e_2)=0.
\]
It compares with the non-trivial connection in Proposition
\ref{hei}. For $n=2$
\[
\gamma(e_1)=\left(
\begin{matrix}
0 & 0 &  0 & 0\\
0 & 0 &  0 & 0\\
a & b & 0 & 0 \\
b & c & 0 & 0
\end{matrix}
\right), \quad \gamma(e_2)=\left(
\begin{matrix}
0 & 0 &  0 & 0\\
0 & 0 &  0 & 0\\
b & c & 0 & 0 \\
c & d & 0 & 0
\end{matrix}
\right), \quad \gamma(e_3)=0=\gamma(e_4),
\]
where all
coefficients in the matrices above are real numbers.

The resulting Lie algebra $\mathfrak h=\mathbb R^{2n} \ltimes_\gamma
\mathbb R^{2n}$ is a 2-step nilpotent Lie algebra. In fact, take
$e_1, e_2, \hdots, e_{2n}, v_1, v_2, \hdots, v_{2n}$ as a basis of
$\mathfrak h$. Then the non-trivial Lie bracket relations are given
by
\begin{equation}\label{rcll}
 {[e_i, v_j]}=\sum_{k=1}^na^{\ k}_{ij}v_{n+k}
 \end{equation}
 for all $1\leq i, j\leq n$,
 where for each $i$ $(a^{\ k}_{ij})$ is the matrix $A_i$.
As the commutator of $\mathfrak h$ is generated by  $\{v_{n+1},
v_{n+2}, \hdots, v_{2n}\}$, which belongs to the center of $\mathfrak h$,
$\mathfrak h$ is 2-step nilpotent. It is
 non-abelian if at least one of the matrices $A_i$ is not trivial.
By Theorem \ref{tower}, $\mathfrak h$ carries a hypersymplectic
structure.

If all the coefficients $a^{\ k}_{ij}$ are rational numbers, the
simply connected Lie group $H$ with Lie algebra $\mathfrak h$ admits
a co-compact lattice. The resulting space may be considered as a
generalization of the Kodaira-Thurston surface, as seen from the
perspective of a hypersymplectic manifold in \cite{Kamada}.

\

 {\bf Acknowledgment\ }\  G. P. Ovando is partially supported by CONICET,
Secyt UNC, ANPCyT and DAAD. \ Y. S. Poon is partially supported by
NSF DMS-0906264.

\

\end{document}